\documentclass[12pt]{amsart}
\usepackage[T1]{fontenc}
\usepackage{lmodern}
\usepackage[margin=1in]{geometry}
\usepackage{amsmath,amssymb,amsthm,mathtools}
\usepackage{microtype}
\usepackage{enumitem}
\usepackage{xcolor}
\usepackage[colorlinks=true,linkcolor=blue!40!black,
  citecolor=blue!40!black,urlcolor=blue!40!black]{hyperref}
\hypersetup{
  pdftitle={The squared antipode of a twisted function algebra is locally unipotent},
  pdfauthor={Pavel Etingof and Shlomo Gelaki},
  pdfsubject={Tensor functors, cotriangular Hopf algebras, and Newton sums}
}
\numberwithin{equation}{section}
\newtheorem{theorem}{Theorem}[section]
\newtheorem{lemma}[theorem]{Lemma}
\newtheorem{corollary}[theorem]{Corollary}
\newtheorem{proposition}[theorem]{Proposition}
\theoremstyle{remark}
\newtheorem{remark}[theorem]{Remark}
\newcommand{\C}{\mathbb C}
\newcommand{\Q}{\mathbb Q}
\newcommand{\Z}{\mathbb Z}
\newcommand{\T}{\mathbb T}
\newcommand{\OO}{\mathcal O}
\newcommand{\Rep}{\operatorname{Rep}}
\newcommand{\Vect}{\operatorname{Vec}_{\mathrm{fd}}}
\newcommand{\Sym}{\operatorname{Sym}}
\newcommand{\tr}{\operatorname{tr}}
\newcommand{\id}{\operatorname{id}}
\newcommand{\ch}{\operatorname{ch}}
\newcommand{\ev}{\operatorname{ev}}

\newcommand{\tors}{\mathrm{tors}}

\setlist[enumerate]{label=\textup{(\roman*)},leftmargin=*,
  itemsep=3pt,topsep=4pt}
\title[The squared antipode of a twisted function algebra]{The squared antipode of a twisted function algebra is locally unipotent}
\author{Pavel Etingof}
\author{Shlomo Gelaki}
\date{}

\begin{document}
\begin{abstract}
Let $G$ be a proalgebraic group over an algebraically closed field $k$
of characteristic zero, and let $F:\Rep(G)\to\Vect(k)$ be a
dimension-preserving tensor functor, possibly nonsymmetric. We prove
that a tensor automorphism $T$ of $F$ is locally unipotent if and only
if its trace on $F(V)$ equals $\dim V$ for every $V\in\Rep(G)$.
Applied to the Drinfeld element of a dual twisted function algebra $(\OO(G)^J)^*$,
this proves the unipotency conjecture of \cite{EG} that the square of
the antipode in $\mathcal{O}(G)^J$ is locally unipotent. Without the trace assumption, we show that 
for any $V\in \Rep(G)$ the
eigenvalues of $T$ on $F(V)$ determine, modulo roots of unity,
its eigenvalues on $F(\mathbb S^\nu V)$ for every Schur functor
$\mathbb S^\nu$. The proof combines
the Koszul identity, polynomial growth of symmetric powers, Fourier
inversion, and congruences for Newton sums.
\end{abstract}
\maketitle

\section{Introduction}\label{sec:introduction}

\subsection{Main results}
Let $k$ be an algebraically closed field of characteristic
zero. For a proalgebraic group $G$ over $k$, let $\Rep(G)$ denote
its category of finite-dimensional rational representations.
Let $\Vect(k)$ denote the category of finite-dimensional $k$-vector
spaces, and let $F:\Rep(G)\to\Vect(k)$ be a tensor functor
\cite[Definition~4.2.5]{EGNO}. The functor $F$ need not be symmetric.
For a tensor automorphism $T$ of $F$, write $T_V:F(V)\to F(V)$ for
the corresponding operator on $F(V)$, $V\in\Rep(G)$.

\begin{theorem}\label{thm:main}
Suppose that $F$ is dimension-preserving. If
\begin{equation}\label{eq:trace}
  \tr(T_V)=\dim V\qquad\text{for every }V\in\Rep(G),
\end{equation}
then $T_V$ is unipotent for every $V$, i.e., $(T_V-1)^{\dim V}=0$.
\end{theorem}

We apply Theorem~\ref{thm:main} to prove the unipotency conjecture
\cite[Conjecture~6.1]{EG}. Let $J$ be a Hopf $2$-cocycle on $\OO(G)$
and equip the Hopf algebra $A=\OO(G)^J$ with the cotriangular form
$R_J=J_{21}^{-1}*J$. Write $S$ for its antipode and $u\in A^*$ for
its Drinfeld character. For a finite-dimensional right $A$-comodule $W$
with coaction $\rho_W:W\to W\otimes A$, set
$u_W=(\id_W\otimes u)\circ\rho_W$. Let $C(W)\subset A$ be the
coefficient coalgebra of $W$.

\begin{corollary}\label{cor:hopf}
Every $u_W$ is unipotent, and $S^2$ is unipotent on every
finite-dimensional subcoalgebra of $A$. More specifically, if
$d=\dim W\geq1$, then
\begin{equation}\label{eq:bounds}
  (u_W-1)^d=0,\qquad
  (S^2-1)^{2d-1}\big|_{C(W)}=0.
\end{equation}
\end{corollary}

Corollary~\ref{cor:hopf} also gives a positive answer to
\cite[Question~7.4]{EG} for twisted function algebras, since all
eigenvalues of $S^2$ are $1$.\footnote{The case of cocycles obtained
by deformation of the trivial cocycle was settled in
\cite[Theorem~6.4]{EG}. For unipotent affine algebraic groups, Brown and Gelaki prove
$S^2=1$ over $\C$ \cite[Corollary~5.6(1)]{BG}. This extends to any $k$ of characteristic zero,
since $G$ and the values of $J$ on a countable basis of
$\OO(G)^{\otimes 2}$ can be defined over a countably generated
subfield of $k$, which embeds in $\C$.}
For general cotriangular Hopf algebras, Bichon gave an example in which
$S^2$ has an eigenvalue of infinite order
\cite[Remark~5]{Bi}.\footnote{For a general cotriangular Hopf algebra,
the dimensions of symmetric powers in its comodule category can grow
exponentially.}

Recall that a cotriangular Hopf algebra $A$ is
\emph{pseudoinvolutive} if
$\tr(S^2|_C)=\dim C$ for every finite-dimensional subcoalgebra
$C\subset A$. Combining Corollary~\ref{cor:hopf} with
\cite[Theorem~3.3]{EG} gives the following characterization.

\begin{corollary}\label{cor:pseudo}
A cotriangular Hopf algebra over $k$ is pseudoinvolutive if
and only if the square of its antipode is locally unipotent.
\end{corollary}

Without the trace assumption, we show that for every Schur functor
$\mathbb S^\nu$ the eigenvalues of $T_{\mathbb S^\nu V}$ agree modulo
torsion with those of $\mathbb S^\nu(T_V)$
(Proposition~\ref{rem:schur}). This means that after identifying
eigenvalues whose ratio is a root of unity and adding their
multiplicities, the two lists agree. Remark~\ref{rem:higher_torsion}
shows why this qualification is needed even for the exterior square.

\subsection{Strategy of proof of Theorem~\ref{thm:main}}
Fix $V$ of dimension $n$ and let $\lambda_1,\ldots,\lambda_n$ be
the eigenvalues of $T_V$, counted with multiplicity. For $1\leq j\leq n$,
consider the two operators
\[
\begin{aligned}
  T_{\wedge^j V}&:F(\wedge^j V)\longrightarrow F(\wedge^j V),\\
  \wedge^j(T_V)&:\wedge^j F(V)\longrightarrow\wedge^j F(V).
\end{aligned}
\]
The first is the component of $T$ at the representation $\wedge^j V$.
The second is the ordinary exterior power of the linear map $T_V$.
The eigenvalues of $\wedge^j(T_V)$ are the products $\lambda_{i_1}\cdots\lambda_{i_j}$
with $i_1<\cdots<i_j$. If the two operators had the same trace for
each $j$, then \eqref{eq:trace} applied to $\wedge^j V$ would give
\[
  \sum_{i_1<\cdots<i_j}\lambda_{i_1}\cdots\lambda_{i_j}
     =\binom nj.
\]
Equivalently, $\prod_{i=1}^n(1+\lambda_i t)=(1+t)^n$, which would
force every $\lambda_i$ to equal $1$. The difficulty is that under
the identification $F(V^{\otimes j})\cong F(V)^{\otimes j}$, the
image of the antisymmetrizer on $V^{\otimes j}$ need not
be the usual antisymmetrizer on $F(V)^{\otimes j}$.

We use the polynomial growth of $\dim F(\Sym^m V)$ to overcome this
difficulty. Let $\Gamma\subset k^\times$ be the group generated by
$\lambda_1,\ldots,\lambda_n$. We record the spectra of the
$T_{\wedge^j V}$ in a polynomial $E(t)$ over $\Z[\Gamma]$.
Each eigenvalue $\gamma$ of $T_{\wedge^j V}$ contributes
$[\gamma]t^j$ once for each occurrence, where $[\gamma]$ is the
basis element indexed by $\gamma$. The Koszul identity identifies
$E(-t)^{-1}$ with the series recording the spectra of the
$T_{\Sym^m V}$. Positivity and polynomial growth imply that, for every
character $\chi:\Gamma\to\C^\times$ with values on the unit circle,
substituting $\chi(\gamma)$ for $[\gamma]$ gives a polynomial all of
whose roots lie on the unit circle. Fourier inversion then bounds the
integral coefficients of the Newton sums of $E$, so that Frobenius
congruences at sufficiently large primes become exact identities.
Together with \eqref{eq:trace}, these yield the rational function
identity \eqref{eq:trace_poles}. Comparing the two sides at their
pole at $1$ forces every $\lambda_i$ to equal $1$.

\subsection{Organization of the paper}
Section~\ref{sec:preliminaries} collects the arithmetic and combinatorial 
ingredients and proves the main technical lemma.
Section~\ref{sec:proofs} proves the main results and records
the consequence for Schur functors.

\subsection*{Acknowledgements}
The main argument of the paper is due to ChatGPT 6 Pro. The authors guided ChatGPT to simplify the argument and then edited and
checked the paper.

\section{Arithmetic and combinatorial preliminaries}
\label{sec:preliminaries}

\subsection{Spectral characters and the Koszul identity}

For an invertible operator $B$ on a finite-dimensional $k$-vector
space $V$ with spectrum in a subgroup $\Gamma\subset k^\times$,
define its \emph{spectral character} by
\[
  \ch_\Gamma(B):=\sum_{\gamma\in\Gamma}m_\gamma[\gamma]
    \in\Z_{\geq0}[\Gamma],
\]
where $m_\gamma$ is the algebraic multiplicity of $\gamma$.
Spectral characters are additive on exact sequences and multiplicative
on tensor products.

Define the augmentation $\varepsilon:\Z[\Gamma]\to\Z$ and the
evaluation $\ev:\Z[\Gamma]\to k$ by
\begin{equation}\label{eq:augev}
  \varepsilon([\gamma])=1,\qquad \ev([\gamma])=\gamma.
\end{equation}
Then
\[
  \varepsilon(\ch_\Gamma(B))=\dim V,\qquad
  \ev(\ch_\Gamma(B))=\tr(B).
\]

Let $F$ and $T$ be as in Theorem~\ref{thm:main}, without assuming
\eqref{eq:trace}. Fix $V$ of dimension $n\geq1$ and let $\Gamma$ be
the subgroup of $k^\times$ generated by the eigenvalues of $T_V$. Since $T$ is a
tensor automorphism, its action on $F(V^{\otimes m})$ identifies with
$T_V^{\otimes m}$. Both exterior and symmetric powers of $V$ are direct
summands of tensor powers, so the spectra of $T_{\wedge^j V}$ and
$T_{\Sym^m V}$ lie in $\Gamma$. Set
\begin{equation}\label{eq:spectralseries}
  \begin{aligned}
  E(t)&:=\sum_{j=0}^n e_jt^j, & e_j&:=\ch_\Gamma(T_{\wedge^j V}),\\
  H(t)&:=\sum_{m\geq0}h_mt^m, & h_m&:=\ch_\Gamma(T_{\Sym^m V}).
  \end{aligned}
\end{equation}
Since $F$ is dimension-preserving, $e_0=[1]$, $e_n$ is a single
basis element, and
\begin{equation}\label{eq:dimensions}
  \varepsilon(e_j)=\binom nj,\qquad
  \varepsilon(h_m)=\binom{n+m-1}{n-1}.
\end{equation}

Exactness of the Koszul complex gives
\[
  \sum_{j=0}^{\min(n,m)}
    (-1)^j[\wedge^jV\otimes\Sym^{m-j}V]=0
  \quad\text{in }K_0(\Rep(G)),\qquad m>0.
\]
Taking spectral characters gives
\begin{equation}\label{eq:koszul}
  E(-t)H(t)=1.
\end{equation}

\subsection{Newton sums and congruences}

Let $\Gamma$ be a finitely generated abelian group. Put
$R:=\Z[\Gamma]$, and for $a\geq1$ define the Adams operation $\psi_a:R\to R$ by
$\psi_a([\gamma])=[\gamma^a]$.
For $n\geq1$ and $E(t)=1+e_1t+\cdots+e_nt^n\in R[t]$,
define its Newton (power) sums $p_m\in R$ by
\begin{equation}\label{eq:newton}
  \frac{tE'(-t)}{E(-t)}=\sum_{m\geq1}p_mt^m.
\end{equation}
Thus $p_1=e_1$, and in any field specialization in which
$E(t)=\prod_{i=1}^n(1+\alpha_it)$, one has
$p_m=\sum_i\alpha_i^m$.

Let $\T:=\{z\in\C:|z|=1\}$. Each character $\chi:\Gamma\to\T$
extends linearly to $\Z[\Gamma]$ and acts coefficientwise on $E(t)$.

\begin{lemma}\label{lem:newton}
Suppose that for every character
$\chi:\Gamma\to\T$, the complex polynomial
$\chi(E(t))$ has degree $n$ and all its roots belong to $\T$.
Let $\Gamma_{\tors}$ be the torsion subgroup of $\Gamma$ and put $D:=(2n)!\,|\Gamma_{\tors}|$. Then
\begin{equation}\label{eq:newton_adams}
  p_{ma}=\psi_a(p_m)\qquad(m,a\geq1,\ \gcd(a,D)=1).
\end{equation}
\end{lemma}

\begin{proof}
Write $p_m=\sum_\gamma c_{m,\gamma}[\gamma]$.
Let $\nu$ be the Haar probability measure on the compact character group
$\widehat\Gamma=\operatorname{Hom}(\Gamma,\T)$. For each $\chi$, the factors
$1+\alpha_i t$ of $\chi(E(t))$ have $|\alpha_i|=1$. Hence
$|\chi(p_m)|=|\sum_i\alpha_i^m|\leq n$, and Fourier inversion gives
\begin{equation}\label{eq:newton_bound}
  |c_{m,\gamma}|=
  \left|\int_{\widehat\Gamma}\chi(p_m)
       \overline{\chi(\gamma)}\,d\nu(\chi)\right|\leq n.
\end{equation}
For every prime $\ell$, the universal Newton identities give
\begin{equation}\label{eq:newton_congruence}
  p_{m\ell}\equiv p_m^\ell\equiv\psi_\ell(p_m)\pmod{\ell R}.
\end{equation}
The first congruence is the identity
$(\sum_i\alpha_i^m)^\ell=\sum_i\alpha_i^{m\ell}$ in characteristic
$\ell$, expressed in elementary symmetric functions. The second is
the Frobenius identity in $\mathbb F_\ell[\Gamma]$.
If $\ell\nmid D$, then $\ell>2n$ and $\gamma\mapsto\gamma^\ell$
is injective on $\Gamma$. Thus $\psi_\ell$ sends distinct basis
elements to distinct basis elements. By \eqref{eq:newton_bound}, every
coefficient of $p_{m\ell}-\psi_\ell(p_m)$ has absolute value at most
$2n<\ell$. By \eqref{eq:newton_congruence}, it is divisible by
$\ell$, so it vanishes.
Factoring $a$ into primes proves \eqref{eq:newton_adams}.
\end{proof}

\subsection{A spectral lemma}

For a finitely generated subgroup $\Gamma\subset k^\times$, put
$\overline\Gamma:=\Gamma/\Gamma_{\tors}$. We also use a bar for the
induced map $\Z[\Gamma]\to\Z[\overline\Gamma]$.

\begin{lemma}\label{lem:spectral}
Let $n\geq1$ and
$E(t)=\sum_{j=0}^n e_jt^j\in\Z_{\geq0}[\Gamma][t]$, where
$e_0=[1]$, $e_n$ is a single basis element, and
$e_1=\sum_{i=1}^n[\gamma_i]$.
Suppose
\begin{equation}\label{eq:growth}
  E(-t)^{-1}=\sum_{m\geq0}h_mt^m,\qquad
  h_m\in\Z_{\geq0}[\Gamma],\qquad
  \limsup_{m\to\infty}\varepsilon(h_m)^{1/m}\leq1.
\end{equation}
Then $\varepsilon(e_j)=\binom nj$, $\varepsilon(h_m)=\binom{n+m-1}{n-1}$,
and the following hold.
\begin{enumerate}
\item\label{item:free}
The image of $E$ modulo torsion is
\begin{equation}\label{eq:free}
  \overline E(t)=
    \prod_{i=1}^n(1+[\overline\gamma_i]t).
\end{equation}
\item\label{item:trace}
If $\ev(E(t))=(1+t)^n$, then
\begin{equation}\label{eq:rigidity}
  e_j=\binom nj[1],\qquad
  h_m=\binom{n+m-1}{n-1}[1].
\end{equation}
\end{enumerate}
\end{lemma}

\begin{proof}
For every character $\chi:\Gamma\to\T$, positivity gives
$|\chi(h_m)|\leq\varepsilon(h_m)$. By \eqref{eq:growth}, the series
$\sum_m\chi(h_m)t^m$ converges for $|t|<1$. Its product with
$\chi(E(-t))$ is $1$, so $\chi(E(t))$ has no root in the open unit
disk. Its constant coefficient is $1$ and its leading coefficient
has modulus $1$. The product of the moduli of its roots is therefore
$1$, so every root lies on the unit circle.

In particular, write $\varepsilon(E(t))=\prod_{i=1}^n(1+\alpha_it)$
with $|\alpha_i|=1$. The linear coefficient is $\varepsilon(e_1)=n$,
so $\sum_i\alpha_i=n$. Taking real parts shows that $\alpha_i=1$ for
every $i$. Thus $\varepsilon(E(t))=(1+t)^n$, and
$\sum_m\varepsilon(h_m)t^m=(1-t)^{-n}$. This proves the two
augmentation formulas.

With $D$ as in Lemma~\ref{lem:newton}, take $m=1$ and $a=1+Dr$ in
\eqref{eq:newton_adams}. Since $p_1=e_1$, we get
\begin{equation}\label{eq:newton_progression}
  p_{1+Dr}=\psi_{1+Dr}(e_1)
    =\sum_{i=1}^n[\gamma_i^{1+Dr}]\qquad(r\geq0).
\end{equation}

For \ref{item:free}, choose a character
$v:\overline\Gamma\to\mathbb R_{>0}$ and put
$v_i:=v(\overline\gamma_i)>0$.
Since $e_n$ is a single basis element, $v(\overline E(t))$
has degree $n$. Write
\[
  v(\overline E(t))=\prod_{j=1}^n(1+\alpha_jt),
  \qquad \alpha_j\ne0.
\]
Applying $v$ after projection to $\overline\Gamma$ in
\eqref{eq:newton_progression} gives
$\sum_j\alpha_j^{1+Dr}=\sum_i v_i^{1+Dr}$ for every $r\geq0$.
Multiplying by $z^r$ and summing over $r$ gives
\begin{equation}\label{eq:free_poles}
  \sum_{j=1}^n\frac{\alpha_j}{1-\alpha_j^Dz}
    =\sum_{i=1}^n\frac{v_i}{1-v_i^Dz}.
\end{equation}
Let $b>0$ occur $m_b$ times among the $v_i$.
Multiply \eqref{eq:free_poles} by $1-b^Dz$ and evaluate at
$z=b^{-D}$ to obtain
\[
  \sum_{\alpha_j^D=b^D}\alpha_j=m_b b.
\]
Each summand has absolute value $b$, so there are at least
$m_b$ such summands. The index sets for distinct $b$ are disjoint, and
$\sum_b m_b=n$, so exactly $m_b$ summands occur for each $b$.
Their sum is $m_b b$, so each has real part $b$ and hence equals $b$.
The multisets $\{\alpha_j\}$ and $\{v_i\}$ therefore coincide.
This proves \eqref{eq:free} after every positive real specialization.
Since these specializations are Zariski dense in the complex
torus associated to $\overline\Gamma$, \eqref{eq:free} holds
in its Laurent polynomial ring.

For \ref{item:trace}, embed the finitely generated field $\Q(\Gamma)$ in
$\C$, so that we may regard the elements of $\Gamma$ as complex numbers. The assumption
$\ev(E(t))=(1+t)^n$ gives $\ev(p_m)=n$ for every $m\geq1$. Thus
\eqref{eq:newton_progression} gives $\sum_i\gamma_i^{1+Dr}=n$ for
every $r\geq0$. Summing as above, we obtain
\begin{equation}\label{eq:trace_poles}
  \sum_{i=1}^n\frac{\gamma_i}{1-\gamma_i^Dz}
    =\frac{n}{1-z}.
\end{equation}
Multiplying by $1-z$ and evaluating at $z=1$ gives
\[
  \sum_{\gamma_i^D=1}\gamma_i=n.
\]
The sum contains at most $n$ terms, each of real part at most $1$.
It follows that all $n$ terms occur and each equals $1$. Thus
$e_1=n[1]$.
For any unitary $\chi$, write
$\chi(E(t))=\prod_{i=1}^n(1+\alpha_it)$, where $|\alpha_i|=1$.
Its linear coefficient is $n$, so $\sum_i\alpha_i=n$ forces
every $\alpha_i=1$. Thus $\chi(E(t))=(1+t)^n$ for every
$\chi$. Fourier inversion gives $e_j=\binom nj[1]$, and the
formula for $h_m$ follows by inverting $E(-t)$.
\end{proof}

\section{Proofs of the main results}\label{sec:proofs}

\begin{proof}[Proof of Theorem~\ref{thm:main}]
Fix $V$ of dimension $n\geq1$ and form $E(t)$ and $H(t)$ as in
\eqref{eq:spectralseries}. By \eqref{eq:koszul}, $H(t)=E(-t)^{-1}$.
Its coefficients are spectral characters, and their augmentations
$\binom{n+m-1}{n-1}$ grow polynomially in $m$ by
\eqref{eq:dimensions}. Thus Lemma~\ref{lem:spectral} applies.
For $0\leq j\leq n$, the trace assumption \eqref{eq:trace} gives
\[
  \ev(e_j)=\tr(T_{\wedge^j V})=\dim\wedge^j V=\binom nj.
\]
Hence $\ev(E(t))=(1+t)^n$. Taking $j=1$ in \eqref{eq:rigidity} gives
$\ch_\Gamma(T_V)=n[1]$. Thus all eigenvalues of $T_V$ are $1$, as
claimed.
\end{proof}

\begin{proof}[Proof of Corollary~\ref{cor:hopf}]
Twisting by $J$ gives a dimension-preserving symmetric tensor
equivalence
$\widetilde F:\Rep(G)\to\operatorname{Comod}_{\mathrm{fd}}(A,R_J)$.
Compose it with the forgetful functor to vector spaces to obtain $F$.
Since $u$ is an algebra character of $A$, the operators $u_W$ form a
tensor automorphism of the forgetful functor. Hence
$X\mapsto u_{\widetilde F(X)}$ defines a tensor automorphism of $F$.
By \cite[equation~(8)]{EG}, the trace of $u_{\widetilde F(X)}$ is the
categorical dimension of $\widetilde F(X)$. The equivalence
$\widetilde F$ preserves this dimension, which in $\Rep(G)$ is
$\dim X$. Theorem~\ref{thm:main} therefore gives the first equality
in \eqref{eq:bounds} for every $W$.

To prove the second equality, write $U=u_W$ and consider the
surjective matrix coefficient map
\[
  c:W^*\otimes W\longrightarrow C(W),\qquad
  c(f\otimes w)=(f\otimes\id_A)\rho_W(w).
\]
By the formula $S^2=u*\id_A*u^{-1}$ of \cite[equation~(6)]{EG},
\begin{equation}\label{eq:matrix-coefficients}
  S^2 c(f\otimes w)=c\bigl((f\circ U)\otimes U^{-1}w\bigr).
\end{equation}
Thus $S^2$ on $C(W)$ is induced by $P\otimes Q$ on $W^*\otimes W$,
where $P=U^*$ and $Q=U^{-1}$. Since $(U-1)^d=0$, we have
$(P-1)^d=0$ and $(Q-1)^d=0$. Now
\[
  P\otimes Q-1=(P-1)\otimes Q+1\otimes(Q-1).
\]
The two summands commute and each has $d$-th power zero, so the
binomial formula gives $(P\otimes Q-1)^{2d-1}=0$. This proves
\eqref{eq:bounds}. Finally, every finite-dimensional subcoalgebra
$C\subset A$ is the coefficient coalgebra of its regular right
comodule. Indeed, the coefficients lie in $C$, and taking $f$ to be
the counit on $C$ gives $c(f\otimes w)=w$ for every $w\in C$.
\end{proof}

\begin{proof}[Proof of Corollary~\ref{cor:pseudo}]
By \cite[equation~(6)]{EG}, $S^2$ preserves every subcoalgebra.
Local unipotence therefore implies the defining trace identities.
Conversely, \cite[Theorem~3.3]{EG} identifies the underlying Hopf
algebra of a pseudoinvolutive cotriangular Hopf algebra with
$\OO(G)^J$. The possible parity modification affects only the
cotriangular form, so Corollary~\ref{cor:hopf} applies to its antipode.
\end{proof}

\subsection{Spectral characters for an arbitrary tensor automorphism}

Let $F,T,V,\Gamma$ be as in Section~\ref{sec:preliminaries}, without
the trace assumption, and let $\lambda_1,\ldots,\lambda_n$ be the
eigenvalues of $T_V$. Write $\ch_{\overline\Gamma}(B)$ for the image
of $\ch_\Gamma(B)$ under $[\gamma]\mapsto[\overline\gamma]$.

\begin{proposition}\label{rem:schur}
For every partition $\nu$, let $\mathbb S^\nu$ and $s_\nu$ denote the
corresponding Schur functor and Schur polynomial. Then
$$\ch_{\overline\Gamma}(T_{\mathbb S^\nu V})
    =s_\nu([\overline\lambda_1],\ldots,[\overline\lambda_n]).$$
\end{proposition}

\begin{proof}
Part~\ref{item:free} of Lemma~\ref{lem:spectral} gives
\begin{equation}\label{eq:exteriorcharacters}
  \sum_{j=0}^n\ch_{\overline\Gamma}(T_{\wedge^j V})t^j
     =\prod_{i=1}^n(1+[\overline\lambda_i]t).
\end{equation}
The Koszul identity \eqref{eq:koszul} now gives
$\overline H(t)=\prod_{i=1}^n(1-[\overline\lambda_i]t)^{-1}$.
Thus $\overline h_m$ is the complete symmetric polynomial of degree
$m$ in the $[\overline\lambda_i]$. Applying spectral characters
modulo torsion to the Jacobi--Trudi identity in $K_0(\Rep(G))$
therefore gives
$\ch_{\overline\Gamma}(T_{\mathbb S^\nu V})
=\det(\overline h_{\nu_i-i+j})_{1\leq i,j\leq\ell(\nu)}
=s_\nu([\overline\lambda_1],\ldots,[\overline\lambda_n])$,
where $\overline h_0=1$ and $\overline h_m=0$ for $m<0$.
\end{proof}

\begin{remark}
\label{rem:higher_torsion}
For exterior squares, the equality in Proposition~\ref{rem:schur} can
fail in $\Z[\Gamma]$. Let $N\geq2$ be even,
choose a primitive $2N$-th root of unity $\xi\in k$, and put $\zeta:=\xi^2$.
On $X:=k^N$, with basis indexed by $a\in\Z/N\Z$, define
$Re_a:=\xi^{-1}e_{a+1},\ 
  Se_a:=\zeta^ae_a,\ Z:=\zeta I$.
Let $G_N:=\langle R,S,Z\rangle\subset\operatorname{GL}(X)$.
The relations $SR=RSZ$, $R^N=Z^{N/2}$, and $S^N=Z^N=1$
show that $G_N$ is finite. Its elements have the unique normal form
$Z^cR^aS^b$, $0\leq a,b,c<N$.
For $V\in\Rep(G_N)$, write $V=\bigoplus_{a,b\in\Z/N\Z}V_{a,b}$,
where $S$ and $Z$ act on $V_{a,b}$ by $\zeta^a$ and $\zeta^b$.
Then $R(V_{a,b})=V_{a+b,b}$.
Equip the ordinary forgetful functor $F$ with the tensorator
$
  j_{V,W}(v\otimes w)=\zeta^{bc}v\otimes w,
  \ v\in V_{a,b},\ w\in W_{c,d}.
$
This defines a dimension-preserving tensor functor, since
$((a,b),(c,d))\mapsto\zeta^{bc}$ is a normalized $2$-cocycle.
Define $T_V(v)=\xi^{b^2}Rv$ on $V_{a,b}$.
The scalar $\xi^{b^2}$ is well-defined modulo $N$ because $N$ is
even, and the identity
$\xi^{(b+d)^2}=\xi^{b^2}\xi^{d^2}\zeta^{bd}$
gives $T_{V\otimes W}j_{V,W}=j_{V,W}(T_V\otimes T_W)$.
Thus $T$ is a tensor automorphism, with
$T_V^N|_{V_{a,b}}=(-1)^{b^2+b}I=I$.

On $X$, the operator $U:=T_X=\xi R$ is the cyclic shift,
so its spectrum is $\mu_N:=\langle\zeta\rangle$, with each
eigenvalue occurring once. Since $Z$ acts on $\wedge^2X$ by
$\zeta^2$, on the underlying vector spaces we have
$ T_{\wedge^2X}=\xi^4\wedge^2R
     =\zeta\,\wedge^2U.$
The multiplicity of $\zeta^j$ in $\wedge^2U$ is the number of
unordered pairs of distinct residues with sum $j$ modulo $N$. Among
the $N$ ordered pairs with this sum, two are diagonal when $j$ is even
and none when $j$ is odd. Thus the multiplicity is $N/2-1$ for even
$j$ and $N/2$ for odd $j$. Multiplication by $\zeta$ interchanges
these multiplicities. Consequently,
\begin{equation}\label{eq:higher_torsion_difference}
  \Delta:=\ch_{\mu_N}(T_{\wedge^2X})
     -\ch_{\mu_N}(\wedge^2U)
     =\sum_{j=0}^{N-1}(-1)^j[\zeta^j]\ne0.
\end{equation}
For a subgroup $H\subseteq\mu_N$ of order $h$, put $q:=N/h$.
Its image in $\Z[\mu_N/H]$ is
\[
  \overline\Delta=
  \begin{cases}
    \displaystyle h\sum_{j=0}^{q-1}(-1)^j
          [\overline\zeta^{\,j}],&q\text{ even},\\[4pt]
    0,&q\text{ odd}.
  \end{cases}
\]
Indeed, its coefficient at $[\overline\zeta^{\,j}]$ is
$\sum_{\ell=0}^{h-1}(-1)^{j+\ell q}$.
Thus, writing $N=2^a m$ with $m$ odd, the discrepancy vanishes
modulo $H$ exactly when $\mu_{2^a}\subseteq H$.
In particular, when $N$ is a power of $2$, it survives quotienting
by every proper subgroup of $\mu_N$.
\end{remark}

\end{document}